\documentclass[12pt]{amsart}

\usepackage{amssymb}

\newtheorem{theorem}{Theorem}[section]
\newtheorem{lemma}[theorem]{Lemma}
\newtheorem{e-proposition}[theorem]{Proposition}
\newtheorem{corollary}[theorem]{Corollary}
\newtheorem{e-definition}[theorem]{Definition\rm}
\newtheorem{remark}{\it Remark\/}

\renewenvironment{proof}{\par\noindent{\bf Proof.}\hspace{0.5em}}
    {\hfill\qed\medskip}
\newenvironment{proofidea}{\noindent{\bf Sketch of proof.}\hspace{0.5em}}
    {\hfill\qed\medskip}

\begin{document}
\title{A New Characterisation of Idempotent States on Finite and Compact Quantum Groups}

\author[UF]{Uwe Franz}
\email{uwe.franz@univ-fcomte.fr}
\author[AS]{Adam Skalski}
\email{a.skalski@lancaster.ac.uk}

\address[UF]{D\'epartement de math\'ematiques de Besan\c{c}on,
Universit\'e de Franche-Comt\'e,
16, route de Gray,
F-25 030 Besan\c{c}on, France, {\rm http://www-math.univ-fcomte.fr/pp\underline{ }Annu/UFRANZ/}}
\address[AS]{Department of Mathematics and Statistics, Lancaster
  University, Lancaster, United Kingdom, {\rm http://www.maths.lancs.ac.uk/\~{ }skalski/}}
\address[AS]{\emph{Permanent address:}
  Department of Mathematics, University of \L\'{o}d\'{z}, ul. Banacha 22,
  90-238 \L\'{o}d\'{z}, Poland}
\thanks{U.F.\ was supported by a Marie Curie Outgoing International
Fellowship of the EU (Contract Q-MALL MOIF-CT-2006-022137) and an ANR Project
(Number ANR-06-BLAN-0015).}

\begin{abstract}
We show that idempotent states on finite quantum groups correspond to
pre-subgroups in the sense of Baaj, Blanchard, and Skandalis.
It follows that the lattices formed by the
idempotent states on a finite quantum group and by its coidalgebras are
isomorphic. We show furthermore that these lattices are also isomorphic for compact quantum
groups, if one restricts to expected coidalgebras.
\end{abstract}

\maketitle

\markboth{Characterisation of Idempotent States}{Franz and Skalski}

\section{Introduction}
\label{sec-intro}

The idempotent measures on a locally compact group are exactly the Haar
measures of its compact subgroups, cf.\ \cite{kawada+ito40,heyer77}. In 1996,
Pal \cite{pal96} has shown that the analogous statement for quantum groups is
false. In \cite{franz+skalski08a}, we have given more examples of
idempotent states on quantum groups that do not come from compact
subgroups. We also provided characterisations of idempotent states on finite
quantum groups in terms of group-like projections
\cite{landstad+vandaele07} and quantum subhypergroups. Subsequently with
Tomatsu we extended some of these
results to compact quantum groups, and determined all idempotent states on the
compact quantum groups $U_q(2)$, $SU_q(2)$, and $SO_q(3)$, cf.\  \cite{franz+skalski+tomatsu09}.

In this note we give a new characterisation of idempotent states on finite quantum groups in terms of the pre-subgroups
introduced in \cite{baaj+blanchard+skandalis99}. That pre-subgroups give rise to idempotent states was not emphasized in
\cite{baaj+blanchard+skandalis99}, but can easily be seen from \cite[Proposition 3.5(a)]{baaj+blanchard+skandalis99}. Here we
prove that, conversely, every idempotent state comes from a pre-subgroup, cf.\ Theorem \ref{main-theo}. As a consequence, we get
a one-to-one correspondence between the idempotent states on a finite quantum group $(\mathsf{A},\Delta)$ and the coidalgebras in
$(\mathsf{A},\Delta)$, cf.\ Corollary \ref{cor-psg-coid}. The isomorphisms providing this bijection have natural explicit
descriptions, cf.\ Remark \ref{rem} after Corollary \ref{cor-psg-coid}. The idempotent states coming from quantum
subgroups are exactly those corresponding to subgroups in the sense of Baaj, Blanchard, and Skandalis, and to coidalgebras of
quotient type, see Proposition \ref{prop-haar}.

The one-to-one correspondence between idempotent states and coidalgebras
extends to compact
quantum groups, if one requires the coidalgebras to be expected, cf.\ Theorem \ref{theo-comp}.

\section{Preliminaries}
\label{sec-prelim}

Recall that a \emph{compact quantum group} is a pair $(\mathsf{A},\Delta)$ of a
unital C${}^*$-algebra $\mathsf{A}$ and a unital $*$-homomorphism
$\Delta:\mathsf{A}\to\mathsf{A}\otimes\mathsf{A}$ such that
$({\rm id}\otimes \Delta)\circ\Delta = (\Delta\otimes{\rm
  id})\circ \Delta$
holds, and the subspaces
\[
{\rm span}\,\{\Delta(b)(\mathbf{1}\otimes a); a,b\in\mathsf{A}\}
\quad\mbox{ and }\quad
{\rm span}\,\{\Delta(b)(a\otimes \mathbf{1}); a,b\in\mathsf{A}\}
\]
are dense in $\mathsf{A}\otimes\mathsf{A}$, cf.\ \cite{woronowicz87,woronowicz98} (here $\otimes$ denotes the minimal tensor
product of $C^*$-algebras reducing to the algebraic tensor product in the finite-dimensional situation). If $\mathsf{A}$ is
finite-dimensional, then $(\mathsf{A},\Delta)$ is called a \emph{finite quantum group} and it admits a counit, i.e.\ a character
$\varepsilon:\mathsf{A}\to \mathbb{C}$ such that $(\varepsilon \otimes {\rm id})\circ \Delta = ({\rm id} \otimes \varepsilon)\circ \Delta
= {\rm id}$. Woronowicz showed that there exists a unique state $h:\mathsf{A}\to\mathbb{C}$ such that
\[
({\rm id}_\mathsf{A}\otimes h)\circ\Delta(a) = h(a)\mathbf{1} =
(h\otimes {\rm id}_\mathsf{A})\circ\Delta(a)\qquad\mbox{ for all } a\in\mathsf{A},
\]
called the Haar state of $(\mathsf{A},\Delta)$. If $(\mathsf{A},\Delta)$ is a finite quantum group, then $h$ is a faithful trace.
A finite quantum group has a unique \emph{Haar element}, i.e.\ a projection $\eta$ such that $\eta a=a\eta=\varepsilon(a)\eta$ for
all $a\in\mathsf{A}$. For more information on finite-dimensional $*$-Hopf algebras and their Haar states, see \cite{vandaele97}.

Define $V:\mathsf{A}\otimes\mathsf{A}\to\mathsf{A}\otimes\mathsf{A}$ by
\[
V(a\otimes b) = \Delta(a)(\mathbf{1}\otimes b)
\]
for $a,b\in\mathsf{A}$. Then $V$
extends to a unitary operator $V:H\otimes H\to H\otimes H$ ($H=L^2(\mathsf{A},h)$ denotes the GNS Hilbert space of the Haar
state), which satisfies $V_{12} V_{13}V_{23} = V_{23}V_{12}$, on $H\otimes H\otimes H$, where we used the leg notation
$V_{12}=V\otimes {\rm
  id}$, etc.
The operator $V$ is called the \emph{multiplicative unitary} of
$(\mathsf{A},\Delta)$, and plays a central role in the approach to quantum
groups developed by Baaj and Skandalis, cf.\ \cite{baaj+skandalis93}.

The notion of a quantum subgroup was introduced by Kac  \cite{kac68} in the
setting of finite ring groups and by Podle\'s \cite{podles95} for matrix pseudo-groups.

\begin{e-definition}\label{def-qsg}
Let $(\mathsf{A},\Delta_\mathsf{A})$ and $(\mathsf{B},\Delta_\mathsf{B})$ be two compact quantum
groups. Then $(\mathsf{B},\Delta_\mathsf{B})$ is called a \emph{quantum subgroup} of
$(\mathsf{A},\Delta_\mathsf{A})$, if there is exists a surjective $*$-algebra
homomorphism $\pi:\mathsf{A}\to\mathsf{B}$ such that
$\Delta_\mathsf{B}\circ \pi = (\pi\otimes \pi)\circ \Delta_\mathsf{A}$.
\end{e-definition}

This definition is motivated by the properties of the restriction map $C(G)\ni f\mapsto
f|_{H}\in C(H)$ induced by a subgroup $H\subseteq G$. If $\mathsf{A}=C(G)$ is
a commutative compact quantum group, then Definition \ref{def-qsg} is equivalent to the usual notion of a closed subgroup.

\begin{e-definition}\label{def-psg}
(\cite[Definition 3.4]{baaj+blanchard+skandalis99}) Let $(\mathsf{A},\Delta_\mathsf{A})$ be a finite quantum group with
multiplicative unitary $V:H\otimes H\to H\otimes H$. Then a \emph{pre-subgroup} of  $(\mathsf{A},\Delta_\mathsf{A})$ is a unit
vector $f\in H$ such that $\varepsilon(f)>0$, and $V(f\otimes f)=f\otimes f$.
\end{e-definition}
Denote by $\mathbf{1}_h\in H$ the cyclic vector that implements the Haar state. For a finite quantum group, $\mathsf{A}\ni
a\mapsto a\mathbf{1}_h\in H$ is an isomorphism and $\varepsilon(f)$ is to be understood via this identification.

We will frequently use this identification and omit $\mathbf{1}_h$ in the rest
of the paper.

A pre-subgroup $f$ is called a \emph{subgroup}, if it belongs to the center of $\mathsf{A}$. In that case $f$ gives rise to a
quantum subgroup in the sense of Definition \ref{def-qsg}, cf.\ Lemma \ref{lem-sg}.

A non-zero element $p\in\mathsf{A}$ in a compact quantum group
$(\mathsf{A},\Delta)$ is called a \emph{group-like projection} \cite[Definition
1.1]{landstad+vandaele07}, if it is a projection, i.e.\ $p^2=p=p^*$, and
satisfies $\Delta(p)(\mathbf{1}\otimes p)=p\otimes p$. We shall see that for finite quantum groups
pre-subgroups and group-like idempotents are essentially the same objects,
i.e.\ that after a rescaling pre-subgroups are group-like projections in
$\mathsf{A}$, cf.\ Corollary \ref{cor-psg-grp}.

For commutative finite quantum groups of the form $\mathsf{A}=C(G)$, pre-subgroups are multiples of indicator
functions of subgroups, cf.\ \cite[Proposition
1.4]{landstad+vandaele07}, but for noncommutative finite quantum groups this
notion is more general than Definition \ref{def-qsg}.

Baaj, Blanchard, and Skandalis defined an order of pre-subgroups by $g\prec f$
if and only if $V (f\otimes g) =f\otimes g$, cf.\ \cite[Proposition 3.7]{baaj+blanchard+skandalis99}.

\section{Characterisations of idempotents states on finite quantum groups}
\label{sec-char}

The coproduct $\Delta:\mathsf{A}\to\mathsf{A}\otimes\mathsf{A}$ leads to
an associative product
\[
\psi_1\star\psi_2=(\psi_1\otimes\psi_2)\circ\Delta
\]
called the
\emph{convolution product}, for linear functionals $\psi_1,\psi_2:\mathsf{A}\to\mathbb{C}$. A state
$\phi:\mathsf{A}\to\mathbb{C}$ is \emph{idempotent}, if
$\phi\star\phi=\phi$. Examples are given by $\phi=h_\mathsf{B}\circ \pi$, if
$(\mathsf{B},\Delta_\mathsf{B})$ is a quantum subgroup of $(\mathsf{A},\Delta_\mathsf{A})$ with
morphism $\pi:\mathsf{A}\to\mathsf{B}$ and Haar state
$h_\mathsf{B}:\mathsf{B}\to\mathbb{C}$. We will call an idempotent state $\phi$ on a compact quantum group
$(\mathsf{A},\Delta)$ a \emph{Haar idempotent state}, if it is of this form.

The natural order
for projections can be used to equip the set of idempotent states on a compact
quantum group with a partial order, i.e.\ $\phi_1\prec\phi_2$ if and only if
$\phi_1\star \phi_2=\phi_2$, cf.\ \cite[Section 5]{franz+skalski08a}

Before we can state and prove the main theorem, we need the following lemma,
which is a variation of \cite[Lemma 4.3]{maes+vandaele98}.

\begin{lemma}\label{lem-grp}
Let $(\mathsf{A},\Delta)$ be a compact quantum group with two states $f$ and
$g$ such that $g\star f=f\star g=f$. Denote by $g_b$ the
functional defined by $g_b(a)=g(ab)$ for $a,b\in \mathsf{A}$. Then we have
\[
f \star g_b = g(b) f
\]
for all $b\in \mathsf{A}$.
\end{lemma}
\begin{proof}
Let $v\in\mathsf{A}$, and set $y=L_f(v)=(f\otimes {\rm id}_\mathsf{A})\circ\Delta(v)$. Then $L_g(y)=(g\otimes{\rm id}_\mathsf{A})\circ\Delta(y)=L_{f\star g}(v) = L_f(v)=y$, therefore
\begin{gather*}
(g\otimes {\rm id}_\mathsf{A})\Big(\big(\Delta(y) - \mathbf{1}\otimes y\big)\big(\Delta(y) - \mathbf{1}\otimes y\big)^*\Big) \\
= L_g(yy^*)-yL_g(y^*) - L_g(y)y^*+yy^*
=  L_g(y^*y)-y^*y,
\end{gather*}
and
\[
(g \otimes f)\Big(\big(\Delta(y - \mathbf{1}\otimes y\big)\big(\Delta(y) - \mathbf{1}\otimes y\big)^*\Big) = 0.
\]
By Cauchy-Schwarz
\begin{gather*}
\left|(g \otimes f)\Big(\big(\Delta(y) - \mathbf{1}\otimes y\big)(b\otimes u)\Big)\right|^2\\
\le
(g \otimes f)\Big(\big(\Delta(y) - \mathbf{1}\otimes y\big)\big(\Delta(y) - \mathbf{1}\otimes y\big)^*\Big)(g\otimes f)(b^*b\otimes u^*u) =0,
\end{gather*}
i.e.\
\[
(g \otimes f)(b\otimes yu)=(g \otimes f)\big(\Delta(y)(b\otimes u)\big)
\]
for all $u,b\in\mathsf{A}$. Recalling the definition of $y$, we get
\begin{gather*}
g(b)(f\otimes f)\big(\Delta(v)(\mathbf{1}\otimes u)\big) \\
= (f\otimes g\otimes f)\Big(\big((\Delta\otimes{\rm id}_\mathsf{A})\circ\Delta(v)\big)(\mathbf{1}\otimes b\otimes u)\Big) \\
= \big((f\star g_b)\otimes f\big)\big(\Delta(v)(\mathbf{1}\otimes u)\big)
\end{gather*}
for all $u,v,b\in\mathsf{A}$. Since ${\rm span}\,\{\Delta(v)(\mathbf{1}\otimes u); u,v\in\mathsf{A}\}$ is dense in
$\mathsf{A}\otimes\mathsf{A}$ and $f$ is nonzero, we get $g(b)f=f\star g_b$.
\end{proof}

For $u,v\in L^2(\mathsf{A},h)$, denote by
$\omega_{u,v}:\mathsf{A}\to\mathbb{C}$ the linear functional $\mathsf{A}\ni
a\mapsto \omega_{u,v}(a) =\langle u, av\rangle = h(u^*av)$.

We have the following characterisation of idempotent states in terms of
pre-subgroups.

\begin{theorem}\label{main-theo}
Let $(\mathsf{A},\Delta)$ be a finite quantum group. Then the map $f\mapsto
\omega_{f,f}$ 
 defines an
order-preserving bijection between the pre-subgroups of $(\mathsf{A},\Delta)$
and the idempotent states on $(\mathsf{A},\Delta)$.
\end{theorem}
\begin{proof}
Let $\omega_{f,f}$ be the state associated to a pre-subgroup $f\in\mathsf{A}$.
We have
\begin{eqnarray*}
(\omega_{f,f}\star \omega_{f,f})(a) &=& \langle f\otimes f, \Delta(a)(f\otimes
f)\rangle \\
&=& \langle f\otimes f, V(a\otimes \mathbf{1}) V^* (f\otimes f)\rangle
\\
&=& \langle f\otimes f,
(a\otimes \mathbf{1}) (f\otimes f)\rangle = \omega_{f,f}(a),
\end{eqnarray*}
for all $a\in \mathsf{A}$, i.e.\ $\omega_{f,f}$ is an idempotent state. This
also follows from \cite[Proposition 3.5(a)]{baaj+blanchard+skandalis99}.

Conversely, let $\phi:\mathsf{A}\to\mathbb{C}$ be an idempotent state. Since
the Haar state is tracial, there exists a unique positive element $\rho_\phi\in\mathsf{A}$ such that
$\phi(a)=\langle\rho_\phi,a\rangle$ for all $a\in\mathsf{A}$. Set
$f_\phi=\sqrt{\rho_\phi}$. Then have $\phi(a)=\langle f_\phi,af_\phi\rangle$ for all
$a\in\mathsf{A}$, and $f_\phi=\sqrt{\rho_\phi}$ is the unique positive element with this property.

By Lemma \ref{lem-grp}, we have $\phi\star\phi_b=\phi(b)\phi$, i.e.\
\begin{eqnarray*}
\langle \rho_\phi\otimes\rho_\phi,a\otimes b\rangle &=& \phi(a)\phi(b) \\
&=& (\phi\star\phi_b)(a) \\
&=&\langle \rho_\phi\otimes
\rho_\phi,\Delta(a)(\mathbf{1}\otimes b)\rangle
\\
&=& \langle \rho_\phi\otimes \rho_\phi,V
(a\otimes \mathbf{1})V^*(\mathbf{1}\otimes b)\rangle \\
&=& \langle V^*(\rho_\phi\otimes \rho_\phi),a\otimes b\rangle
\end{eqnarray*}
for all $a,b\in\mathsf{A}$, since $V(\mathbf{1}\otimes b)=\Delta(\mathbf{1})(\mathbf{1}\otimes b)=\mathbf{1}\otimes b$. Therefore
we have $V(\rho_\phi\otimes \rho_\phi)=\rho_\phi\otimes \rho_\phi$. Recalling the definition of $V$ and the identification
between $H$ and $\mathsf{A}$, this means $\Delta(\rho_\phi)(\mathbf{1}\otimes\rho_\phi)=\rho_\phi\otimes\rho_\phi$. Applying
$\varepsilon$ to the left-hand-side, we get $\rho_\phi^2=\varepsilon(\rho_\phi)\rho_\phi$. Therefore $\varepsilon(\rho_\phi)>0$ and
$f_\phi=\sqrt{\rho_\phi}=\frac{\rho_\phi}{\sqrt{\varepsilon(\rho_\phi)}}$. Clearly, $f_\phi$ is a unit vector,
$\varepsilon(f_\phi)=\sqrt{\varepsilon(\rho_\phi)}>0$, and $V (f_\phi\otimes f_\phi)=f_\phi\otimes f_\phi$, i.e.\ $f_\phi$ is a
pre-subgroup.

Let $g$ be another pre-subgroup with $\phi=\omega_{g,g}$. If we can show $g\ge0$, then this implies $g=f_\phi$. Applying
$\varepsilon$ to $\Delta(g)(\mathbf{1}\otimes g)=g\otimes g$, we get $g^2=\varepsilon(g)g$. Applying $\phi$ to the Haar element $\eta$,
we see $\varepsilon(g)=\varepsilon(f_\phi)$. Furthermore, $\omega_{g,g}=\omega_{f_\phi,f_\phi}$ is equivalent to $gg^*=f_\phi
f_\phi^*$. Therefore we get $||g||=||f_\phi||$, and $g/\varepsilon(g)$ is an idempotent with norm one. This implies that $g$ is an
orthogonal projection, therefore positive, and we see that $f\mapsto \omega_{f,f}$ defines indeed a bijection between the set of
pre-subgroups in $\mathsf{A}$ and the set of idempotent states on $\mathsf{A}$.

Let now $f,g$ be two pre-subgroups such that $g\prec f$, i.e.\ $V (f\otimes g)=f\otimes g$. Then
\begin{eqnarray*}
(\omega_{f,f}\star\omega_{g,g})(a) &=& \langle f\otimes g,\Delta(a)(f\otimes g)\rangle \\
&=& \langle f\otimes g,
V(a\otimes\mathbf{1})V^*(f\otimes g)\rangle \\
&=& \langle f\otimes g, (a\otimes\mathbf{1})(f\otimes g)\rangle \\
&=& \omega_{f,f}(a)
\end{eqnarray*}
for all $a\in\mathsf{A}$, i.e.\ $\omega_{g,g}\prec\omega_{f,f}$.

Conversely, if $\omega_{f,f}\star\omega_{g,g}=\omega_{f,f}$, then $\omega_{g,g}\star\omega_{f,f}=\omega_{f,f}$ since idempotent
states are invariant under the antipode, see \cite[Lemma 5.2]{franz+skalski08a}. By  Lemma \ref{lem-grp} we get
$\omega_{f,f}\star(\omega_{g,g})_b=\omega_{g,g}(b)\omega_{f,f}$, and (recalling that $\omega_{f,f}(\cdot) = \varepsilon(f) \langle
f, \cdot\rangle$ and similarly for $g$)
\begin{eqnarray*}
\langle f\otimes g, a\otimes b\rangle &=& \frac{\omega_{f,f}(a)\omega_{g,g}(b)}{\varepsilon(f)\varepsilon(g)}
= \frac{\big(\omega_{f,f}\star(\omega_{g,g})_b\big)(a)}{\varepsilon(f)\varepsilon(g)} \\
&=&\langle f\otimes g,\Delta(a)(\mathbf{1}\otimes b)\rangle
= \langle V^*(f\otimes g),a\otimes b\rangle,
\end{eqnarray*}
for all $a,b\in\mathsf{A}$, i.e.\ $g\prec f$.
\end{proof}

Note that we have also shown in this proof that any pre-subgroup is
self-adjoint and becomes a projection after an appropriate scaling.

\begin{corollary}\label{cor-psg-grp}
Let $(\mathsf{A},\Delta)$ be a finite quantum group. The map $f\mapsto \frac{f}{\varepsilon(f)}$ defines a bijection between the
pre-subgroups and the group-like projections of $(\mathsf{A},\Delta)$.
\end{corollary}

A \emph{right coidalgebra} $\mathsf{C}$ in a compact quantum group is a unital
$*$-subalgebra $\mathsf{C}\subseteq\mathsf{A}$ such that
$\Delta(\mathsf{C})\subseteq \mathsf{A}\otimes\mathsf{C}$. Baaj, Blanchard,
and Skandalis have shown that the lattice of pre-subgroups of a
finite quantum groups is isomorphic to its lattice of right coidalgebras, cf.\
\cite[Proposition 4.3]{baaj+blanchard+skandalis99}.

\begin{corollary}\label{cor-psg-coid}
Let $(\mathsf{A},\Delta)$ be a finite quantum group. Then the lattice of idempotent states on $(\mathsf{A},\Delta)$ and the
lattice of right coidalgebras in $(\mathsf{A},\Delta)$ are isomorphic.
\end{corollary}

\begin{remark}\label{rem}
We can also give an explicit description of this bijection. Let $\phi:\mathsf{A}\to\mathbb{C}$ be an idempotent state. The one
can show that $T_\phi:\mathsf{A}\to\mathsf{A}$, $T_\phi=({\rm id}_\mathsf{A}\otimes \phi)\circ\Delta$ defines a conditional
expectation, i.e.\ a projection $E:\mathsf{A}\to\mathsf{C}$ onto a unital $*$-subalgebra $\mathsf{C}\subseteq\mathsf{A}$ such
that $||E||=1$, $E(\mathbf{1})=\mathbf{1}$, and $h\circ E=h$. Furthermore, since $T_\phi$ is right-invariant,
$T_\phi(\mathsf{A})$ is a coidalgebra. Conversely, to recover an idempotent state $\phi$ from a right coidalgebra
$\mathsf{C}\subseteq\mathsf{A}$, set $\phi=\varepsilon\circ E_{\mathsf{C}}$, where $E_{\mathsf{C}}$ denotes the unique
$h$-preserving conditional expectation onto $\mathsf{C}$. See also Theorem \ref{theo-comp}.
\end{remark}

\begin{lemma}\label{lem-sg}
Let $(\mathsf{A},\Delta)$ be a finite quantum group, $f$ a subgroup of $(\mathsf{A},\Delta)$, i.e.\ a pre-subgroup that belongs
to the center of $\mathsf{A}$, and put $\tilde{f}=\frac{f}{\varepsilon(f)}$. Then $(\mathsf{A}_f,\Delta_f)$ is a quantum subgroup of
$(\mathsf{A},\Delta)$, with $\mathsf{A}_f=\mathsf{A}f=\{af;a\in\mathsf{A}\}$, and
$\Delta_f:\mathsf{A}_f\to\mathsf{A}_f\otimes\mathsf{A}_f$ and $\pi_f:\mathsf{A}\to\mathsf{A}_f$ given by
\[
\Delta_f(a) = \Delta(a)(\tilde{f}\otimes \tilde{f}) \qquad \mbox{ and } \qquad \pi(a)=a\tilde{f}
\]
for $a\in\mathsf{A}$.
\end{lemma}
\begin{proof}
This follows from Corollary \ref{cor-psg-grp} and
 \cite[Proposition 2.1]{landstad+vandaele07}.
\end{proof}

For any quantum subgroup $(\mathsf{B}, \Delta_\mathsf{B})$ of $(\mathsf{A},\Delta)$,
$\mathsf{A}/\!/\mathsf{B}=\{a\in\mathsf{A}; ((\pi\otimes {\rm id})\circ
\Delta_\mathsf{A})(a)=\mathbf{1}_\mathsf{B}\otimes a\}$ defines a right coidalgebra. A right
coidalgebra is said to be of \emph{quotient type}, if it is of this form.

Using the previous Lemma, one can check that under the one-to-one
correspondences given in Theorem \ref{main-theo} and Corollary
\ref{cor-psg-coid}, Haar idempotent states correspond to subgroups and
coidalgebras of quotient type.
\begin{e-proposition}\label{prop-haar}
Let $\phi$ be an idempotent state. Then the following are equivalent.
\begin{enumerate}
\item
The state $\phi$ is a Haar idempotent state.
\item
The pre-subgroup $f_\phi$ is a subgroup.
\item
The coidalgebra $\mathsf{C}_\phi$ is of quotient type.
\end{enumerate}
\end{e-proposition}

\section{Extension to compact quantum groups}
\label{sec-ext}

For a compact quantum group $(\mathsf{A},\Delta)$, in general the Haar state $h$ is no longer a trace, and for a closed unital
$*$-subalgebra $\mathsf{B}\subseteq\mathsf{A}$ there might exist no $h$-preserving conditional expectation
$E_\mathsf{B}:\mathsf{A}\to\mathsf{B}$. It turns out that the existence of such a conditional expectation is the condition we
have to add to extend the bijection between idempotent states and right coidalgebras. Recall that a compact quantum group is
called coamenable if its reduced version is isomorphic to the universal one (equivalently, the Haar state $h$ is faithful and
$\mathsf{A}$ admits a character, cf.\ \cite[Corollary 2.9]{bedos+murphy+tuset01}). In particular every coamenable compact quantum
group admits a bounded counit.

\begin{theorem}\label{theo-comp}
Let $(\mathsf{A},\Delta)$ be a coamenable compact quantum group. Then there
exists an order-preserving bijection between the expected right coidalgebras in
$(\mathsf{A},\Delta)$ and the idempotent states on
$(\mathsf{A},\Delta)$.
\end{theorem}
\begin{proofidea}
Given an idempotent state $\phi \in \mathsf{A}^*$ we define a completely positive idempotent projection $E_{\phi}=
(\textup{id}_{\mathsf{A}} \otimes \phi)\circ \Delta$. An application of Lemma 3.1 shows that $E_{\phi} (E_{\phi} (a) E_{\phi}
(b)) = E_{\phi} (a) E_{\phi} (b)$ for all $a,b \in \mathcal{A}$, where $\mathcal{A}$ is the $*$-Hopf algebra spanned by
coefficients of the unitary corepresentations of $\mathsf{A}$. Density of $\mathcal{A}$ in $\mathsf{A}$ and the continuity
argument implies that $E_{\phi}(\mathsf{A})$ is an algebra; the right invariance of $E_{\phi}$ expressed by the equality $\Delta
\circ E_{\phi} = (\textup{id}_{\mathsf{A}} \otimes E_{\phi}) \circ \Delta $ implies that $E_{\phi}(\mathsf{A})$ is a right
coidalgebra.

Conversely, if $\mathsf{C}$ is an expected right coidalgebra, let $E_{\mathsf{C}}$ denote the corresponding conditional
expectation. We can show that if $\mathsf{C}'=\{b \in\mathsf{A}: E_{\mathsf{C}}(b) = 0\}$, then for all $\omega \in \mathsf{A}^*,  b
\in \mathsf{C}'$,  $(\omega \otimes \textup{id}_{\mathsf{A}} ) (\Delta(b)) \in
\mathsf{C}'$.
This implies that $E_{\mathsf{C}}$
is right invariant and thus $E_{\mathsf{C}}= (\textup{id}_{\mathsf{A}} \otimes \phi)\circ \Delta$ for the idempotent state
$\phi:= \varepsilon \circ E_{\mathsf{C}}$.
\end{proofidea}

\section*{Acknowledgements}
This work was started while U.F.\ was visiting the Graduate School of
Information Sciences of Tohoku University as Marie-Curie fellow. He would like
to thank Professors Nobuaki Obata, Fumio Hiai, and the other members of the
GSIS for their hospitality. We would also like to thank Eric Ricard and Reiji
Tomatsu for helpful comments and suggestions.


\end{document}